\documentclass[reqno]{amsart}
\usepackage{amsfonts,latexsym,amscd}
\usepackage{amsmath,amsthm}
\usepackage{amssymb,easymat}
\newtheorem*{theorem}{\sc Theorem} 
\newtheorem*{proposition}{\sc Proposition}

\theoremstyle{definition}

\newtheorem*{remark}{\sc Remark}

\newcommand{\p}{\omega}
\begin{document}

\title{Connections on Metriplectic Manifolds}
\author{Daniel Fish}
\address{Department of Mathematics \& Statistics, Portland State University,
Portland, OR, U.S.} \email{djf@pdx.edu}

\begin{abstract} 
In this note we discuss conditions under which a linear connection on a manifold equipped with both a symmetric (Riemannian) and a skew-symmetric (almost-symplectic or Poisson) tensor field will preserve both structures.
\end{abstract}
\maketitle

 If $(M,g)$ is a (pseudo-)Riemannian manifold, then classical
results due to T. Levi-Civita, H. Weyl and E. Cartan \cite{Schouten} show that
for any $(1,2)$ tensor field $T^{i}_{jk}$ which is skew-symmetric by lower indices, there exists a unique linear connection $\Gamma$ preserving the metric ($\nabla^{\Gamma} g=0$), with $T$ as its torsion tensor: $T^{i}_{kj}=\frac{1}{2}(\Gamma^{i}_{jk}-\Gamma^{i}_{kj})$. 
 
 It has also been shown \cite{Gelfand} that given any symmetric (by lower indices) $(1,2)$ tensor $S^{i}_{jk}$ on a symplectic manifold $(M,\omega)$, there exists a unique linear connection preserving $\omega$ which has $S$ as its symmetric part, i.e., $S^{i}_{jk} = \frac{1}{2}(\Gamma^{i}_{jk}+\Gamma^{i}_{kj})$. Moreover, it is known  \cite{Vaisman} that if $\p$ is a regular Poisson tensor on $M$, then there always exists a linear connection on $M$ with respect to which $\p$ is covariantly constant. Such connections are called \emph{Poisson} connections, and can be chosen to coincide with the Levi-Civita connection of the metric $g$ (if $g$ is Riemannian) in certain cases.\par

Considering these results, one is naturally led to the question: Given a skew-symmetric $(0,2)$ tensor $\p$,  and a (pseudo-)Riemannian metric $g$ on a manifold $M$, when do there exist linear connections preserving $\p$ and $g$ simultaneously: 
 \begin{equation}\label{parallel}
\nabla^{\Gamma} \p + \nabla^{\Gamma} g = 0\;?
\end{equation}

Motivated by the terminology of P.J. Morrison \cite{Mor}, we call the a manifold equipped with both a (pseudo-)Riemannian metric $g$ and a skew-symmetric $(2,0)$ tensor $P$ a \emph{metriplectic} manifold, and a connection which preserves both tensors will be called a \emph{metriplectic connection}. In the first section we restrict ourselves to the case in which both $\p=P^{-1}$ and $g$ are nondegenerate, that is $\p$ is almost-symplectic and $g$ is Riemannian. We combine the results from \cite{Schouten} and \cite{Gelfand} to derive a necessary condition for a connection $\Gamma$ to be a metriplectic connection. We also discuss the form of $\Gamma$ in the almost-Hermitian and symplectic cases. The main result of this section is the following

\noindent \textbf{Proposition} \emph{Any connection $\Gamma$ with symmetric part $\Pi$ and torsion $T$ that preserves both a Riemannian metric $g$ and an (almost-)symplectic form $\p$ has the form $\Gamma = \Pi + T$
\[\Pi = L(g) + \bar{g}(T),\,\,\,\,\,T= \p^{-1}\nabla^{g}\p-(1/2)\p^{-1}d\p  - \bar{\p}\bar{g}(T),\]}where $L(g)$ is the Levi-Civita connection defined by the metric $g$, and the ``bar" operators $\bar{g}(T)$ and $\bar{\p}\bar{g}(T)$ are related to the symmetries of the torsion $T$. 

In the second section we give the proof of a theorem due to Shubin \cite{Shubin} which states that if $M$ admits a metriplectic connection, and $P=\p^{-1}$ is nondegenerate, then $M$ is a K\"{a}hler manifold. We also formulate an observation made by Vaisman \cite{Vaisman} as a generalization of Shubin's theorem in the case that $P$ is degenerate.

\section{\textbf{Necessary Conditions on the Metriplectic Connection $\Gamma$}}
We consider here the case where $\omega=P^{-1}$ is skew-symmetric and nondegenerate (not necessarily closed), and $g$ is Riemannian. Suppose that $\Gamma = \Pi + T$ is a connection on the Poisson-Riemannian manifold $(M,g,\omega)$ with $\Pi $ and $T$ symmetric and skew-symmetric (torsion) tensors respectively. Assume that $\Gamma$ satisfies (\ref{parallel}).  Since $\nabla^{\Gamma} g=0$, we know from \cite{Gelfand} that the symmetric tensor $\Pi$ must have local components
\begin{equation} \label{pi} \Pi^{i}_{jk} = L(g)^{i}_{jk} + g^{is}(T^{q}_{sj}g_{qk} + T^{q}_{sk}g_{jq}),\end{equation}
where $L(g)^{i}_{jk} = (1/2)g^{is}(g_{js,k}+g_{ks,j}-g_{jk,s})$ is the Levi-Civita connection for $g$. On the other hand, since $\nabla^{\Gamma}\omega=0$, we have (\cite{Schouten})
\begin{equation} \label{Tor} T^{i}_{jk} = L(\p)^{i}_{jk} - \p^{is}(\Pi^{q}_{sj}\p_{qk} + \Pi^{q}_{sk}\p_{jq}), \end{equation}
where $L(\p)^{i}_{jk} =  (1/2)\p^{is}(\p_{js,k}+ \p_{sk,j}+ \p_{jk,s} )$. We introduce the following operator: for any nondegenerate $(0,2)$ tensor $h$ define the linear operator $\bar{h}$ on $(1,2)$ tensors by
\[\bar{h}(B)^{i}_{jk} = h^{is}(B^{q}_{sj}h_{qk} + B^{q}_{sk}h_{kq}).\]
With this definition, we can write (\ref{pi}) and (\ref{Tor}) as
\[\Pi = L(g) + \bar{g}(T)\;\;\;\textrm{and}\;\;\;T = L(\p) - \bar{\p}(\Pi).\]
Thus, the original connection $\Gamma$ has the form
\begin{eqnarray}\label{gamma} \nonumber \Gamma &=& L(g)+L(\p)-\bar{\p}(\Pi)+\bar{g}(T),\\
                                         &=&L(g)+L(\p)-\bar{\p}(L(g))-\bar{\p}(\bar{g}(T))+\bar{g}(T).
                                                 \end{eqnarray}
                                                 
Notice that when $\bar{h}$ operates on a connection form $A$, the result is related to the covariant derivative of $h$ with respect to $A$ as follows: 
\begin{equation}\label{bar} \nabla^{A}h = \partial{h} - h (\bar{h}(A))\;\;\textrm{or}\;\;\ \bar{h}(A)=h^{-1}\partial{h}-h^{-1}\nabla^{A}h.\end{equation}
In particular,  $\,\, \bar{\p}(L(g))= \p^{-1}(\partial{\p}-\nabla^{g}\p)$
where $\nabla^{g}$ is the covariant derivative with respect to $g$. So we have 
\[\Gamma = L(g)+L(\p)- \p^{-1}(\partial{\p}-\nabla^{g}\p)- \bar{\p}(\bar{g}(T))+\bar{g}(T).\]
Rewriting $L(\p)$ as $\p^{-1}\partial{\p}-(1/2)\p^{-1}d\p$, we have the following 

\begin{proposition}Any connection $\Gamma$ the preserves both a Riemnannian metric $g$ and an (almost-)symplectic form $\p$ has the form $\Gamma = \Pi +T$ with
\[\Pi = L(g) + \bar{g}(T),\,\,\,\,\,T= \p^{-1}\nabla^{g}\p-(1/2)\p^{-1}d\p  - \bar{\p}\bar{g}(T).\] If $\p$ is closed (i.e. $(M,\p)$ is symplectic),then $\Gamma = L(g) + \bar{g}(T) + \p^{-1}\nabla^{g}\p-\bar{\p}\bar{g}(T)$.\end{proposition}

\subsection{\sc Almost-Hermitian Connections with Totally Skew Torsion} Suppose now that $g$ and $\p$ are related by an almost-complex structure $J = g^{-1}\p$ satisfying $J^{2}=-I$. Observe that $\bar{g}(T)=0$ if and only if $T$ is \emph{totally skew-symmetric} with respect to the metric $g$ (that is, $T_{ijk} = T^{q}_{ij}g_{qk}$ is an exterior 3-form). In this case, \begin{equation}\label{herm}T = \p^{-1}\nabla^{g}\p-(1/2)\p^{-1}d\p.\end{equation}
Thus $\Gamma$ reduces to the canonical almost-hermitian connection with totally skew torsion (see e.g. \cite{Gauduchon} \cite{Lich}).
Indeed, the torsion 3-form $T(X,Y,Z) = <T(X,Y),Z>_{g}$ can be expressed in terms of the Nijenhuis tensor $N(X,Y,Z) = <N(X,Y),Z>_{g}$ of the almost-complex structure $J$ (see \cite{Ivan}) as follows. 
\begin{proposition} If the torsion of an almost-Hermitian connection $\Gamma$ is totally skew-symmetric, then
\[T(X,Y,Z) = (1/2)N(X,Y,Z) - (1/2)d\p(JX,JY,JZ)\]
for all $X,Y,Z$.
\end{proposition}
\begin{proof}
\[2T_{ijk}= (2\p^{qn}\nabla^{g}_{n}\p_{ij} - \p^{qn}d\p_{ijn})g_{qk}\]
\[= 2J^{n}_{k}\nabla^{g}_{n}\p_{ij} - J^{n}_{k}(\nabla^{g}_{n}\p_{ij} + \nabla^{g}_{j}\p_{ni} + \nabla^{g}_{i}\p_{jn})\]
\[= J^{n}_{k}\nabla^{g}_{n}\p_{ij} -J^{n}_{k}\nabla^{g}_{j}\p_{ni} -J^{n}_{k}\nabla^{g}_{i}\p_{jn}\]
\[= J^{n}_{k}\nabla^{g}_{n}\p_{ij} -J^{n}_{i}\nabla^{g}_{j}\p_{kn} +J^{n}_{j}\nabla^{g}_{i}\p_{kn}\]
\[= -J^{n}_{k}\nabla^{g}_{n}\p_{ji}-J^{n}_{i}\nabla^{g}_{n}\p_{kj} -J^{n}_{j}\nabla^{g}_{n}\p_{ik}+ J^{n}_{i}\nabla^{g}_{n}\p_{kj} -J^{n}_{j}\nabla^{g}_{n}\p_{ki} -J^{n}_{i}\nabla^{g}_{j}\p_{kn} +J^{n}_{j}\nabla^{g}_{i}\p_{kn}\]
\[= -J^{n}_{k}\nabla^{g}_{n}\p_{ji}-J^{n}_{i}\nabla^{g}_{n}\p_{kj} -J^{n}_{j}\nabla^{g}_{n}\p_{ik}+ N_{ijk}.\]
Using the fact that $\p_{in}\nabla^{g}_{q}J^{n}_{j} = -J^{n}_{j}\nabla^{g}_{q}\p_{in}$, we have 
\[J^{n}_{k}\nabla^{g}_{n}\p_{ji} = J^{n}_{k}g_{jr}\nabla^{g}_{n}J^{r}_{i} = J^{n}_{k}(J^{s}_{j}\p_{sr})\nabla^{g}_{n}J^{r}_{i} = -J^{n}_{k}J^{s}_{j}J^{r}_{i}\nabla^{g}_{n}\p_{sr} .\]
Permuting the indices $i,j,k$ gives us
\[2T_{ijk} = (\nabla^{g}_{t}\p_{sr} + \nabla^{g}_{r}\p_{ts} + \nabla^{g}_{s}\p_{rt})J^{r}_{i}J^{s}_{j}J^{t}_{k} + N_{ijk} \]
\[= -d\p_{rst}J^{r}_{i}J^{s}_{j}J^{t}_{k} + N_{ijk}.\]
\end{proof}

Using this expression for $T$ together with (\ref{herm}), it is easy to see that $M$ is:
\begin{itemize}
\item[] Hermitian ($J$ is integrable) $\leftrightarrow$ $T(X,Y,Z)=-(1/2)d\p(JX,JY,JZ)$, \item[] Symplectic ($\p$ is closed) $\leftrightarrow$  $T = (1/2)N$, and 
\item[] K\"{a}hler  ($\p$ is closed and $g$-parallel)  $\leftrightarrow$ $\Gamma$ is the Levi-Civita connection determined by $g$. \end{itemize}

\subsection{\sc The Symplectic case}
As mentioned above, if $M,\p$ is symplectic, then the connection $\Gamma$ takes the form 
\begin{equation}\label{symp} \Gamma = L(g) + \bar{g}(T) + \p^{-1}\nabla^{g}\p-\bar{\p}\bar{g}(T).\end{equation}
Applying (\ref{bar}) to $\p$ and $\Gamma$, we see that the condition $\nabla^{\Gamma}\p = 0$ is equivalent to
\begin{eqnarray*}0 & = & \partial{\p} - \p(\bar{\p}(\Gamma))\\
                    &=& \partial{\p} - \p(\bar{\p}(\Pi)) - \p(\bar{\p}(T))\\
                   &=& \partial{\p} - \p(\bar{\p}(L(g))) - \p(\bar{\p}(\bar{g}(T)))- \p(\bar{\p}(T))\\ 
                    &=& \nabla^{g}\p - \p(\bar{\p}(\bar{g}(T)))- \p(\bar{\p}(T)).
\end{eqnarray*}
However,  $0 = \nabla^{g}\p -\p(\bar{\p}(\bar{g}(T)))-\p(T)$. Thus we obtain the following condition on $T$:
\[\bar{\p}(T)=T \;\;\;\textrm{or}\;\;\;(\bar{\p} - I)T = 0.\]
\begin{remark}[] \emph{This result can be derived directly from the condition $\nabla^{\Gamma}\p=0$ together with the Jacobi condition for $\p$, and is easily seen to be equivalent to the following cyclic condition on the indices of the tensor $T\p$, \[T^{s}_{ij}\p_{sk}+T^{s}_{ki}\p_{sj}+T^{s}_{jk}\p_{si}=0.\]}\end{remark}
\noindent Thus, we may substitute $\bar{g}(\bar{\p}(T))$ for $\bar{g}(T)$ in (\ref{symp}). Writing the difference $\bar{g}(\bar{\p}(T))-\bar{\p}(\bar{g}(T))$ as $[\bar{g},\bar{\p}](T)$, we arrive at the following
\begin{proposition} A linear connection \ $\Gamma$ on a symplectic-Riemannian manifold $(M,g,\p)$ which preserves both $g$ and $\p$ has the form 
\[\Gamma = L(g) + \p^{-1}\nabla^{g}\p+ [\bar{g},\bar{\p}](T).\]
\end{proposition}

\section{\textbf{Necessary Conditions on the structure of $(M,g,P)$}}

It is known \cite{Lich} that if $M$ is a K\"{a}hler manifold, then the K\"{a}hler form is parallel with respect to the Levi-Civita connection $L$ on $M$ defined by the K\"{a}hler metric, in which case (\ref{parallel}) clearly holds (with $\Gamma = L$). In this section we will discuss a partial converse to this fact due to M. Shubin.
\subsection{\sc{Shubin's Theorem}}
On a manifold $M$, let $g_{0}$ be a Riemannian metric and let $\p=P^{-1}$ be an almost-symplectic form (non-degenerate and skew-symmetric). Let $L(g_{0})$ denote the Levi-Civita connection associated with $g_{0}$, and let $K = g_{0}+\p$. We denote the covariant derivatives with respect and $L(g_{0})$ by $\nabla^{0}$. The following theorem is a reformulation of a result by Shubin \cite{Shubin}, and its proof follows Shubin's proof, with some variations. 

\begin{theorem}
If \ $(M,g_{0},\p)$ is an almost-symplectic Riemannian manifold, and $\nabla^{0} K = 0$, then there exists a complex structure $J$ on $M$ such that the metric $g$ defined by \ $g(X,Y) = \p(X,JY)$ is parallel with respect to $g_{0}$ (thus $L(g) = L(g_{0})$), and defines a K\"{a}hler structure on $M$. 
\end{theorem}
\begin{proof}
First note that if  $\nabla^{0}K = 0$, then $\nabla^{0}\p =0$. Since $\nabla^{0}$ is symmetric, it follows (see remark 1.4 in \cite{Gelfand}) that $d\p = 0$, and so $\p$ is symplectic. Now, any symplectic manifold \cite{Shubin} admits an almost-complex structure $J$ defined by $J=A(-A^{2})^{-1/2}$, where $A$ is the linear operator defined by
\[g_{0}(AX,Y) = \p(X,Y).\]
Since both $g_{0}$ and $\p$ are parallel with respect to $\nabla^{0}$, it is clear that the operator $A$ will also be parallel, thus $\nabla^{0} J=0$. The integrability of $J$ then follows from the expression
\[N_{J}(X,Y) = (\nabla^{0}_{JX}J)Y - (\nabla^{0}_{JY}J)X + J(\nabla^{0}_{Y}J)X - J(\nabla^{0}_{X}J)Y\]
for the Nijenhuis torsion of $J$ (see \cite{DeLeon}).

The metric $g(X,Y) = \p(X,JY)$ is Hermitian with respect to $J$. Therefore, it defines a K\"{a}hler structure on $M$. Furthermore, the equality $g_{ij}=\p_{ik} J^{k}_{j}$ shows that $\nabla^{0}g=0$. The connection $L(g_{0})$ is symmetric and compatible with $g$, so it must coincide with the Levi-Civita connection $L(g)$.
\end{proof}

\subsection{\sc Generalization to a degenerate $P$}
If the tensor $P$ is degenerate, then we cannot construct the covector $\p = P^{-1}$ on $M$. In order to deal with this possibility, we change setting from the cotangent to the tangent bundle.

With $g_{0}$ and $\nabla^{0}$ as above, suppose that $M$ is equipped with a (possibly degenerate) Poisson tensor $P$, and let $K=g_{0}^{-1} + P$. If $\,\nabla^{0} K=0$, then $\nabla^{0} P=0$ and $M$ is a regular Poisson manifold with symplectic foliation $\mathcal{S}(M)$ defined by the kernel of $P$ (see \cite{Vaisman}). The restriction $P_{S}$ of $P$ to a symplectic leaf $S$ is nondegenerate, and $S$ is endowed with a symplectic form $\p = P_{S}^{-1}$.

A classical result of Lichnerowicz \cite{Lich} states that there exist local coordinates $x^{i}$ along $\mathcal{S}(M)$ and $y^{i}$ along $\mathcal{N}$ ( the transverse foliation orthogonal to $\mathcal{S}(M)$) in which $g_{0}$ and $\p$ have the form $g_{0} = g' + g''$ where
\[ g' = (g_{0})_{ij}(y)dy^{i}dy^{j}, \;\;\; g'' = (g_{0})_{ij}(x)dx^{i}dx^{j},\;\;\;  \p = \p_{ij}(x)dx^{i}\wedge dx^{j}. \]

By restricting $\p$ and $g''$ to a symplectic leaf $S$, we are in the situation described by Shubin's Theorem above. Thus, we have a complex structure $J$ which defines a Hermitian metric $g_{s} = g_{ij}(x)dx^{i}dx^{j}$ on $S$ which is parallel with respect to $g''$ (and, therefore, with respect to $g_{0}$). We can extend $J$ \ by $0$ to all of $M$, and define a new metric $\tilde{g}$ on $M$ by the formula $\tilde{g} = g' + g_{s}$.

This metric is called a \emph{partially K\"{a}hler} metric. It is parallel with respect to $g_{0}$ and, when restricted to the symplectic leaf $S$, is a Hermitian metric on $S$. In his book \cite{Vaisman}, Vaisman concludes from these remarks that ``the parallel Poisson structures of a Riemannian manifold $(M,g_{0})$ are exactly those defined by the K\"{a}hler foliation of the $g_{0}$-parallel partially-K\"{a}hler metrics of $M$ (if any)". One can view this statement as the following generalization of Shubin's Theorem.

\begin{theorem} If \ $\nabla^{0} K=0$, then $M$ is a regular Poisson manifold (with the Poisson tensor $P$), and there exists a complex structure $J$ on the symplectic leaves of $M$ such that the metric $\tilde{g}$ defined above is parallel with respect to $g_{0}$ (thus $\nabla^{\tilde{g}} = \nabla^{0}$), and the restriction of \ $\tilde{g}$ to the symplectic leaf \ $S$ defines a K\"{a}hler foliation on $M$.
\end{theorem}

\section{\textbf{Related Questions}}

We have shown that the preservation of both  Riemannian and  Poisson structures on $M$ by a linear connection imposes certain conditions on the connection itself, as well as on the structure of the manifold $(M,g,P)$. In future work we will address some related questions, including: When can one guarantee the existence of a metriplectic connection on a manifold $M$, and is there an optimal or canonical choice of such a connection, similar to the canonical connection given in \cite{Gauduchon}?

\bibliographystyle{plain}
\bibliography{MetriplecticConnections}
\end{document}